\documentclass{article}
\usepackage{a4wide}
\usepackage{amssymb}
\usepackage{amsfonts}
\usepackage{amsxtra}
\usepackage{latexsym}
\usepackage{graphicx}
\usepackage{xcolor} 
\usepackage{wasysym}
\usepackage{authblk}
\usepackage{amsmath}

\usepackage{tikz}

\newtheorem{theo}{Theorem}[section]
\newtheorem{theo*}{Theorem}
\newtheorem{lemma}[theo]{Lemma}

\newtheorem{cor}[theo]{Corollary}
\newtheorem{cor*}{Corollary}

\newcommand{\N}{{\mathbb{N}}}
\newcommand{\Z}{{\mathbb{Z}}}

\newcommand{\R}{{\mathbb{R}}}
\newcommand{\C}{{\mathbb{C}}}

\newcommand{\F}{{\mathbb{F}}}

\newcommand{\qed}{\hspace*{\fill}$\Box$}

\begin{document}

%%%%%%%%%%%%%%%%%%%%%%%%%
% Title.
%%%%%%%%%%%%%%%%%%%%%%%%%

\title{Two-colorings of finite grids: variations on \\ a theorem of Tibor Gallai}

% Authors in alphabetical order
\author{
Bogdan Dumitru\thanks{University of Bucharest, Faculty of Mathematics and Informatics. \texttt{bogdan.dumitru@fmi.unibuc.ro}}
\,
Mihai Prunescu\thanks{Research Center for Logic, Optimization and Security (LOS), Faculty of Mathematics and Computer Science, University of Bucharest, Academiei 14, Bucharest (RO-010014), Romania; \texttt{mihai.prunescu@gmail.com}} \thanks{Simion Stoilow Institute of Mathematics of the Romanian Academy, Research unit 5, P. O. Box 1-764, RO-014700 Bucharest, Romania. \texttt{mihai.prunescu@imar.ro}}
}

\date{}
\maketitle

%%%%%%%%%%%%%%%%%%%%%%%%%
% Abstract.
%%%%%%%%%%%%%%%%%%%%%%%%%

\begin{abstract}
\parindent 0 cm 
A celebrated but non-effective theorem of Tibor Gallai states that for any finite set $A$ of $\Z^n$ and for any finite number of colors $c$ there is a minimal $m$ such that no coloring of the finite $m^n$-grid can avoid that a homothetic image of $A$ is monochromatic. We find (or confirm) $m$ for equilateral triangles, squares, and various types of rectangles. Also, we extend the problem from homothety to general similarity, or to similarity generated using some special rotations. In particular, we compute Gallai similarity numbers for lattice rectangles similar to $1\times k$ (in all orientations) for $k=2,3,4$. The solutions have been found in the framework of the Satisfiability Problem in Propositional Logic (SAT). While some questions were solved using managed brute force, for the more computationally intensive questions we used modern SAT solvers together with symmetry breaking techniques. Some other minor questions are solved for triangles and squares, and new lower bounds are found for regular hexagons on the triangular lattice and for three-dimensional cubes in $\Z^3$. 
\vspace{1mm}


\vspace{1mm}

{\bf M.S.C.-Classification}: 05C15, 05B05, 52C05, 68V15 
\end{abstract}

%%%%%%%%%%%%%%%%%%%%%%%%%%%%%%%%%%%%%%%%%%%%%%%%%%%%%%%%%%%%%%%%%%%%%%%%%%%%%%%%%%%%%%%%%%%%%%%%%%%%%%%%%%%%%%%%%%% 

\section{Introduction} 

According to a Theorem by van der Waerden, for all natural numbers $k$ and $l$ there exists a natural number $n(k,l)$ such that no matter how we paint the natural numbers in  an arbitrary segment of the axis of length $n(k,l)$ with $k$ colors, there exists a monochromatic arithmetic progression of length $l$, see \cite{vdW}. A very nice proof is presented in Khinchin's book \cite{K}. We observe that some arithmetic progression of length $l$ is both similar and homothetic with the standard arithmetic progression $\{0, 1, \dots, l-1\}$, because in dimension $1$ these notions coincide. 

 A set $B \subset \Z^n$ is a $\Z$-{\bf homothetic} image of a set $A \subset \Z^n$ if there exist $y \in \Z^n$ and $x \in \N \setminus \{0\}$ such that $B = xA+y$.  In general, we will use only the term {\bf homothetic}. 

\begin{theo}[Gallai]\label{Gallai}
Let $A$ be a finite subset of $\Z^n$. Then any finite coloring of $\Z^n$ contains a monochromatic subset $B$ which is $\Z$-homothetic with $A$.
\end{theo} 

There are different proofs, like R.~Rado's \cite{R} and a remarkably short one by E.~Witt \cite{W}. See also Soifer's wonderful monograph \cite{S}.  In all proofs, exactly like in the Theorem of van der Waerden, it is shown that there is a natural number $m$ depending on $A$ and on the number of colors $c$ such that no matter how one colors the finite grid $\{0, 1, \dots, m-1\}^n$, this set contains a monochromatic subset $B$ which is homothetic with $A$. 

Unfortunately, the proofs are not computationally effective. From the proof one can derive only huge upper bounds for the corresponding $m$. On the other hand, the complexity of the problem makes the search for the minimal value of $m$ difficult.

Let ${\bf \Gamma (A, c)} \in \N$ be the minimal natural number such that for any $c$-coloring of the finite grid $\{0, 1, \dots, \Gamma(A, c) - 1\}^n$, we cannot avoid that some subset $B$ of this finite grid, homothetic with $A$, is monochromatic. Given the edge-length $m$ of a finite grid, the condition that all subsets homothetic with $A$ aren't monochromatic can be written down as a boolean formula. If $m \geq \Gamma(A, c)$, the propositional formula expressing this condition is {\bf unsatisfiable}, while the propositional formulae corresponding to $m < \Gamma(A, c)$ are {\bf satisfiable}. We call $\Gamma(A, c)$ the {\bf Gallai Homothety Number} for the set $A$ and for $c$ colors. 

From the geometric point of view, not only homotheties are interesting, but also similarities. A {\bf similarity} is a function $\sigma : \R^n \rightarrow \R^n$ such that $\sigma = \tau \circ \delta \circ \rho$, where $\rho$ is a rotation, $\delta$ is a dilation, and $\tau$ is a translation. 

As a homothety is a particular case of similarity, Gallai's Theorem implies that for every finite subset $A \subset \Z^n$ and for every number of colors $c$ there exists a minimal value ${\bf \Sigma(A, c)}$ such that no coloring of the finite grid $\{0, 1, \dots, \Sigma(A, c) - 1\}^n$ can avoid that some finite subset $B$ of this finite grid, which is similar with $A$, is monochromatic. The similar subsets induce more constraints and in general one has:
$$\Sigma(A, c) \leq \Gamma(A,c). $$

We call $\Sigma(A, c)$ the {\bf Gallai Similarity Number} of $A$ for $c$ colors. We address the following questions: given $A$ and $c$, find the Gallai numbers $\Sigma(A, c)$ and $\Gamma(A, c)$. To solve the questions we look for the smallest $m$ for which the corresponding propositional formula is unsatisfiable. For a general account on satisfiability in propositional logic, see the handbook Biere et. al. \cite{Biere_HandbookSAT_2009}.

 We will treat only the case of $2$-colorings, but the situation is similar for any $c \geq 2$. 

In most experiments, the number of dimensions is $n = 2$.
In most cases, $A$ will be a triangle or a square. In order to better visualize the configurations, we will consider equilateral triangles in triangular grids.  

For both problems, we apply {\bf managed brute force}. This method can be described as follows. The grid $G_m$ is embedded in $G_{m+1}$ by adding a row consisting of $m+1$ equidistant points in the triangular case, respectively a new row and a new column consisting of $m+1$ lattice points in the square case - this means $2m+1$ new lattice points. For any solution $s : G_m \rightarrow \{0,1\}$ that avoids some specific monochromatic configurations, we write down all solutions $s' : G_{m+1} \rightarrow \{0,1\}$ extending $s$, which avoid the given monochromatic configurations in the bigger grid - meaning that $s'\,|\,G_m = s$. Of course, only configurations containing the new lattice points have to be checked. If $m + 1 = \Gamma(A, c)$, respectively $m+1 = \Sigma(A, c)$, no such extension does exist anymore. We will denote the finite triangular grid containing $m(m+1)/2$ points with $T_m$, respectively the finite square grid containing $m^2$ points with $S_m$. 

The managed brute force is sufficient for the $2$-color similarity and homothety questions in triangles, and for the $2$-color similarity question in squares. However, the method doesn't work anymore for the homothety question in squares, because of the computation time, as long as one uses only one processor. Possibly, one could accelerate the managed brute force computation by using a computer network and running parallel extension jobs on different processors or different computers. Instead, we used for this question modern SAT-solving with symmetry breaking. So we found out that $\Gamma(\square, 2) = 15$, confirming results found by Walton and Li in \cite{WL} and Bacher and Eliahou \cite{bacher2010extremal}. 

If one tries to estimate $\Gamma(\Box, 2)$ from Witt's proof, one finds only an upper bound in the order of thousands. This shows to what extent such proofs can be ineffective. 

For $k \geq 2$ let $R_k$ be the rectangle $R_k = \{0,1\} \times \{0, k\}$. We find $\Gamma(R_k, 2)$ for $k = 2, 3, 4, 5$. We also address the question of prohibiting homothetic images of both $R_k$ and $R_k' = \{0,k\} \times \{0,1\}$. Moreover, we compute $\Sigma(R_k,2)$ for rectangles similar to $R_k$ (in all orientations) for $k=2,3,4$.
Some other works on rectangles address a different problem. In \cite{FGGP} and \cite{SP} one studies many-colored rectangular grid such that no horizontal rectangle is monochromatic. This approach differs from ours, as we look only for $2$-colorings of square grids avoiding monochromatic homothetic images of only one type of unit-rectangle, in order to complement Gallai's Theorem. 

We also apply the SAT-based approach to regular hexagons in the triangular lattice and to axes-parallel cubes in $\Z^3$, obtaining non-trivial lower bounds for the corresponding minimally unsatisfiable configurations. 

{\bf Open question 1}: {\it Are there other methods to compute the Gallai numbers?}

{\bf Open question 2}: {\it To what extent one can find human-readable proofs for our results?} 

As the Ramsey theory is the branch of combinatorics that studies the appearance of order in large structures, amid  apparent disorder, see \cite{GRS}, looking for the exact Gallai numbers is a problem naturally matching in this framework. Other interesting technique to be applied is given by the pseudo-boolean programming, see \cite{HR}. Indeed, the condition that a square $ABCD$ is not monochrome, can be written as:
$$0 < s(A) + s(B) + s(C) + s(D) < 4,$$
where the additions take place in $\Z$ and the values of $s$ are $0$ or $1$.

%%%%%%%%%%%%%%%%%%%%%%%%%%%%%%%%%%%%%%%%%%%%%%%%%%%%%%%%%%%%%%%%%%%%%%%%%%%%%%%%%%%%%%%%%%%%%%%%%%%%%%%%%%%%%%%%%%%

\section{Linear problems}

We start with an easy question. {\it Is it possible to color the planar square lattice $\Z \times \Z$ in two colors, say white and black, such that every axes-parallel square whose vertices
are lattice points has an odd number of black vertices?} Consider the
square grid $S_3 = \{0, 1, 2\} \times \{0, 1, 2\}$. We associate to the involved lattice points boolean variables
$x_1, x_2, \dots, x_9$, such that $x_i=0$ means that the lattice point in question is colored
white. The corresponding points form a matrix as follows:
\[
\begin{pmatrix}
x_1 & x_2 & x_3 \\
x_4 & x_5 & x_6 \\
x_7 & x_8 & x_9
\end{pmatrix}
\]
The condition that every axes-parallel square has an odd number of black vertices can be translated in the following system of linear equations over the field
$\F_2$:
\begin{eqnarray*}
{\textcolor{blue}{x_1}} + x_2 + x_4 + x_5 &=& 1\\
x_2 + {\textcolor{blue}{x_3}} + x_5 + x_6 &=& 1\\
x_4 + x_5 + {\textcolor{blue}{x_7}} + x_8 &=& 1\\
x_5 + x_6 + x_8 + {\textcolor{blue}{x_9}} &=& 1\\
x_1 + x_3 + x_7 + x_9 &=& 1
\end{eqnarray*}
We add the first, the second, the third and the fourth equation together. It follows that:
\begin{eqnarray*}
x_1 + x_3 + x_7 + x_9 &=& 0,
\end{eqnarray*}
and this contradicts the fifth equation. We proved the following:

\begin{theo}
    $S_2$ is the maximal finite square grid $S_m$ that accepts a $2$-coloring $s : S_m \rightarrow \{0,1\}$ such that all horizontal squares with vertices in lattice points have an odd number of black vertices. 
\end{theo}

Let's try another problem of this type. {\it We look for maximal grids $S_m = \{0, 1, \dots, m-1\} \times \{0, 1, \dots, m-1\}$ and colorings $c : S_m \rightarrow \{0,1\}$ such that every homothetic image of the unit square $\{0,1\} \times \{0,1\}$ has exactly two white and two black vertices.} 

For $m \geq 2$, the set of homothetic images of the unit square, also described as axes-parallel or horizontal squares, form the following set:
\[
  \Box_m = \{\,\,\{(i, j), (i+a, j), (i +a, j + a), (i, j+a)\} \subseteq S_m\,\mid \, i, j, a \in \N, a\neq 0 \}.
\]

As a $2$-coloring of $\Z \times \Z$ is a function $s : \Z \times \Z \rightarrow \{0,1\}$,
we will consider $2$-colorings of the square $m$-grid, $s : S_m \rightarrow \{0,1\}$. In
order to write down boolean constraints for colorings, we associate to $S_m$ the set of
$m^2$ boolean variables $\{s(x,y) \,\mid \, 0 \leq x, y < m\}$.

The formula
\[
  \delta(s(A), s(B), s(C), s(D)) := [ s(A) + s(B) + s(C) + s(D) = 2 ],
\]
where the addition is considered in $\Z$, expresses
the fact that the square $ABCD$ has exactly two vertices of each color. We consider the family of boolean formulas:
\[
  \Pi_m := \bigwedge _{ABCD \in \Box_m} \delta(s(A), s(B), s(C), s(D)),
\] 
$\Pi_m$ is a boolean formula of $m^2$ variables, as $\delta(a, b, c, d)$ can be written also as 
$$\delta(a, b, c, d)=\bar a \bar b cd \vee \bar a b \bar c d \vee \bar a bc \bar d \vee a \bar b \bar c d \vee a \bar b c \bar d \vee a b \bar c \bar d,$$
with conjunction written as juxtaposition and negation of $x$ written as $\bar x$. $\Pi_m$ is satisfiable if and only if there is a coloring $c : S_m \rightarrow \{0, 1\}$ with the given property. The maximal solution is found by managed brute force: 

\begin{theo}\label{th:Pi}
The formulas $\Pi_m$ are satisfiable if and only if $m \leq 4$. Consequently, $S_4$ is the maximal grid accepting a $2$-coloring such that all horizontal squares with vertices on the grid have exactly $2$ black vertices. 
\end{theo}

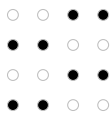
\begin{figure}[ht]
\centering
\begin{tikzpicture}[scale=0.4]
    % 1 = Black (Fill), 0 = White (Draw)
    \foreach \x/\y in {2/3, 3/3, 0/2, 1/2, 2/1, 3/1, 0/0, 1/0} 
       \fill[black] (\x, \y) circle (0.18);
    \foreach \x in {0,...,3}
        \foreach \y in {0,...,3}
           \draw[gray!50] (\x, \y) circle (0.18);
\end{tikzpicture}
\caption{A 2-coloring of the $4 \times 4$ grid such that every horizontal square with vertices on the grid has exactly two white and two black vertices. This is a maximal solution for this problem.}
\label{fig:square4}
\end{figure}

{\bf Proof}: There are $6$ solutions for $\Pi_2$, $8$ solutions for $\Pi_3$, $8$ solutions for $\Pi_4$. See Figure \ref{fig:square4} for a maximal solution. 
$\Pi_5$ is unsatisfiable. \qed  

%%%%%%%%%%%%%%%%%%%%%%%%%%%%%%%%%%%%%%%%%%%%%%%%%%%%%%%%%%%%%%%%%%%%%%%%%%%%%%%%%

\section{Monochromatic subsets}\label{sec:monochromatic}

Let $S$ be a finite set of points and let ${\cal A} $ be a family consisting of some subsets of $S$. {\it We want to find a boolean (propositional) formula expressing the fact that under some coloring $s : S \rightarrow \{0,1\}$, no subset $A \in {\cal A}$ is monochromatic.} To this goal, we introduce boolean variables $\{x_i \mid i \in S\}$ for every element of $S$. To some set $A$ we write down the condition that $A$ is not monochromatic:
$$\varepsilon(A) = \varepsilon^+(A) \wedge \varepsilon^-(A) = \left ( \bigvee_{i \in A} x_i \right ) \wedge \left ( \bigvee_{i \in A} \neg x_i\right ) .$$
Indeed, $\varepsilon^+(A)$ means that at least one $x_i$ is $1$ while $\varepsilon^-(A)$ means that at least one $x_i$ is $0$. Now we define the boolean (propositional) formula:
$$\Upsilon = \Upsilon(S, {\cal A}) = \bigwedge_{A\in {\cal A}} \varepsilon(A). $$
The following statement is evident:
\begin{lemma}\label{lemma:formula}
    The propositional formula $\Upsilon(S, {\cal A})$ is unsatisfiable if and only if for every $2$-coloring $s : S \rightarrow \{0,1\}$ there exists $A \in {\cal A}$ which is monochromatic. 
\end{lemma} 

In the following sections, the set $S=T_m$ or $S = S_m$ will be the finite triangular, square or cubic grid with $m$ rows and ${\cal A} = {\cal A}_m$ will be a family of subsets, like all horizontal equilateral triangles with vertices on the grid, all horizontal or skew squares, and so on. The corresponding formula:
$$\Upsilon_m = \Upsilon(S_m, {\cal A}_m)$$ 
will express the fact that we can avoid monochrome subsets of interest. By Gallai's Theorem, there will be always a minimal $m_0$ such that the corresponding formula $\Upsilon_{m_0}$ is unsatisfiable. All formulas $\Upsilon_k$ with $k \geq m_0$ are unsatisfiable as well. Our goal is to find the corresponding $m_0$ for various families of subsets. We will use the notation $\Upsilon(S, {\cal A})$ as an operator which associates a propositional formula with a pair $(S, {\cal A})$ where $S$ is finite and $\cal A$ is a given family of subsets of $S$. 

When we use SAT-solvers, we take advantage of the fact that $\Upsilon(S, {\cal A})$ is already in conjunctive normal form (CNF). So, in preparing input data, one adds for every $A \in {\cal A}$ the two disjunctive {\bf clauses} $\varepsilon^+(A)$ and $\varepsilon^-(A)$ on the clause list. 

%%%%%%%%%%%%%%%%%%%%%%%%%%%%%%%%%%%%%%%%%%%%%%%%%%%%%%%%%%%%%%%%%%%%%%%%%%%%%%%

\section{Triangles}\label{sec:triangular} 

This section is dedicated to equilateral triangles. They are natural subsets of the following plane triangular lattice.
The triangular lattice is $T = \langle 1, \theta \rangle$, where $\theta = e^{\frac{\pi i}{3}}$. For $m \geq 1$, the triangular $m$-grid is the set
\[
  T_m = \{ a + b \theta \mid a, b \in \N,\, a + b < m\}.
\]
The set $T_m$ has $m(m+1)/2$ elements. In order to address its elements, it is convenient to introduce a level-oriented system of coordinates. According to this system, every
level contains the points $(x,y)$ where $0 \leq y < m$ is the index of the level and
$0 \leq x \leq y$ are the ordinates of the $y+1$ points of the level $y$. The convex closure of $T_m$ is an equilateral triangle with edge-length $m - 1$.

The set $E_m$ of all equilateral triangles with vertices in $T_m$ can be given as follows:
\[
 E_m = \{\,\,\{(i, j), (i - a, j + b), (i + b, j + a + b)\} \subseteq T_m \,\mid \, i, j \in \N,\, a, b \in \Z, (a,b) \neq (0,0) \},
\]
where the vertices are given according to the level-oriented coordinate system. We are also interested in those equilateral triangles whose edges are parallel with the axes and are upwards oriented:
\[
  \Delta_m = \{\,\,\{(i, j), (i+a, j), (i, j - a)\} \subseteq T_m \,\mid \, i, j, a\in \N, a \neq 0\},
\]
and in the equilateral triangles whose edges are parallel with the axes and are downwards oriented:
\[
  \nabla_m = \{\,\,\{(i, j), (i - a, j - a), (i, j - a)\} \subseteq T_m \, \mid \, i, j, a \in \N, a \neq 0\}.
\]
$\Delta_m$ and $\nabla_m$ are subsets of $E_m$. 

We introduce three families of propositional formulas:
\[
  \Phi_m := \Upsilon(T_m, E_m),
\]
\[
  \Gamma_m := \Upsilon(T_m, \Delta_m \cup \nabla_m), 
  \qquad
  \Psi_m := \Upsilon(T_m, \Delta_m) .
\]
The satisfiability of $\Phi_m$ corresponds to the existence of a $2$-coloring of $T_m$ with
no monochromatic equilateral triangle {\it similar} with the unit triangle of the grid. The formulas $\Psi_m$ encode the corresponding restrictions for axes-parallel upward oriented triangles, meaning the {\it homothetic} case. The formulas $\Gamma_m$ say that axes-parallel upward and downward oriented triangles are not monochromatic. 

\begin{figure}[ht]
  \centering
  \begin{tikzpicture}[scale=0.6]
    \foreach \x/\y in {2/3, 1/1.5, 4/0}
      \draw (\x,\y) circle (0.12);
    \foreach \x/\y in {0/0, 2/0, 3/1.5}
      \fill (\x,\y) circle (0.12);
  \end{tikzpicture}
  \caption{This colored grid $T_3$ does not contain any monochromatic horizontal, skew or reversed equilateral triangle. It is a maximal solution for this problem.} \label{fig:triangle3}
\end{figure}
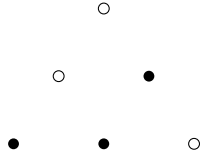

\begin{theo}\label{th:Phi}
The formulas $\Phi_m$ are satisfiable if and only if $m \leq 3$.   Consequently $\Sigma(\triangle, 2) = 4$. One cannot color $T_4$ without the occurrence of a monochromatic horizontal, skew or reversed equilateral triangle. 
\end{theo}

{\bf Proof}: Using managed brute force, we find $6$ solutions of $\Phi_2$ and $18$ solutions for $\Phi_3$. See Figure \ref{fig:triangle3} for a maximal solution. The formula $\Phi_4$ is unsatisfiable. 
\qed 

\begin{theo}\label{th:GammaPsi}
The formulas $\Gamma_m$ and $\Psi_m$ are satisfiable if and only if $m \leq 4$. Consequently $\Gamma(\triangle, 2) = 5 $. One cannot color $T_5$ without the occurrence of a horizontal upwards-oriented monochromatic equilateral triangle. This is the best bound for the homothety class of horizontal upwards oriented equilateral triangles. 
\end{theo}

\begin{figure}[ht]
  \centering
  \begin{tikzpicture}[scale=0.6]
    \foreach \x/\y in {3/4.5, 2/3, 5/1.5, 2/0, 4/0}
      \draw (\x,\y) circle (0.12);
    \foreach \x/\y in {4/3, 1/1.5, 3/1.5, 0/0, 6/0}
      \fill (\x,\y) circle (0.12);
  \end{tikzpicture}
  \caption{This colored grid $T_4$ does not contain any monochromatic horizontal upwards oriented triangle. It is a maximal solution for this problem. We observe that the grid contains a skew equilateral triangle. Of course, as $\Sigma(\triangle, 2) = 4$, one couldn't avoid this. } \label{fig:triangle4}
\end{figure}
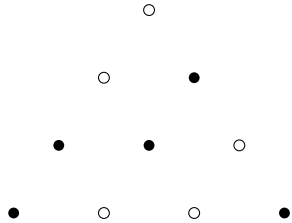

{\bf Proof}:  Using managed brute force for $\Gamma_m$, we find $6$ solutions for $m=2$, $18$ solutions for $m=3$, $36$  solutions for $m=4$. See Figure \ref{fig:triangle4} for a maximal solution. We have the same situation  for $\Psi_m$. Moreover, for $m \leq 4$, $\Gamma_m$ and $\Psi_m$ have the same solutions. $\Gamma_5$ and $\Psi_5$ are unsatisfiable.
\qed 

%%%%%%%%%%%%%%%%%%%%%%%%%%%%%%%%%%%%%%%%%%%%%%%%%%%%%%%%%%%%%%%%%%%%%%%%%%%%%%%%%%%%%%%%%%%%%%%%%%% 

\section{Squares}\label{sec:squares}

The square lattice is $S = \langle 1, i \rangle \subseteq \C$, which we identify with
$\Z \times \Z$. For $m \geq 1$, the square $m$-grid is the set
\[
  S_m = \{a + bi \, \mid \, a, b \in \N,\, a < m,\, b < m\}.
\]
In Cartesian coordinates we have
\[
  S_m = \{0, 1, \dots, m -1\} \times \{0, 1, \dots, m-1\}.
\]
The set $S_m$ has $m^2$ elements, and the convex closure of $S_m$ is a square with edge-length $m-1$.

The set $F_m$ of all squares with vertices in $S_m$ can be defined as follows:
\[
  F_m = \{\,\,\{(i, j), (i + b, j - a), (i + a + b, j + b - a), (i + a, j + b)\}
  \subseteq  S_m \,\mid\, i, j \in \N,\, a, b \in \Z, (a,b)\neq(0,0) \}.
\]
We are also interested in those squares which are parallel with the lattice axes:
\[
  \Box_m = \{\,\,\{(i, j), (i+a, j), (i +a, j + a), (i, j+a) \}\subseteq S_m \,\mid \, i, j, a \in \N, a \neq 0 \}.
\]

We consider two families of boolean formulas:
\[
  \Omega_m := \Upsilon(S_m, F_m), 
\]
\[
  \Xi_m := \Upsilon(S_m, \Box_m).
\]
$\Omega_m$ encode colorings with no monochromatic horizontal or skew square, corresponding to the {\it similarity} question.  $\Xi_m$ encode colorings with no monochromatic horizontal square,  corresponding to the {\it homothety} question, as originally stated in Gallai's Theorem. 

The managed brute force method can be applied to the family $(\Omega_m)$, in analogy with the triangular case. In this case the method leads to a complete classification of the maximal solutions.

\begin{theo}\label{th:Omega}
The formulas $\Omega_m$ are satisfiable if and only if $m \leq 6$. Consequently $\Sigma(\square, 2) = 7$. One cannot color $S_7$ without the occurrence of a horizontal or skew monochromatic square. 
\end{theo}

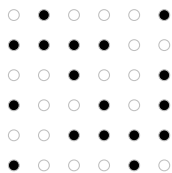
\begin{figure}[ht]
\centering
\begin{tikzpicture}[scale=0.4]
    \foreach \x/\y in {1/5, 5/5,
    0/4, 1/4, 2/4, 3/4,
    2/3, 5/3,
    0/2, 3/2, 5/2,
    2/1, 3/1, 4/1, 5/1,
    0/0, 4/0} 
       \fill[black] (\x, \y) circle (0.18);
    \foreach \x in {0,...,5}
        \foreach \y in {0,...,5}
           \draw[gray!50] (\x, \y) circle (0.18);
\end{tikzpicture}
\caption{A 2-coloring of the $6 \times 6$ grid with no monochromatic horizontal or skew squares. This is a maximal solution for this problem.}
\label{fig:square6}
\end{figure}

{\bf Proof}: There are $14$ solutions for $\Omega_2$, $248$ solutions for $\Omega_3$, $5006$ solutions for $\Omega_4$, $7120$ solutions for $\Omega_5$, $56$ solutions for $\Omega_6$. See Figure \ref{fig:square6} for a maximal solution.  $\Omega_7$ is unsatisfiable. \qed 

For $\Xi_m$ the managed brute force method no longer performs well once $m$ becomes moderate. The authors succeeded in exhibiting solutions for $\Xi_{12}$ but this method alone is not sufficient to decide whether $\Xi_{13}$ and $\Xi_{14}$ are satisfiable. 

As the managed brute force method reached  its limits, we went to more professional SAT-solving.

We employed modern CDCL solvers, specifically \texttt{CaDiCaL} for generating unsatisfiability proofs and \texttt{Kissat} for deep search on satisfiable instances, see the account \cite{BiereEtAl_Kissat_SATComp2020}. For another useful Python toolkit for prototyping with SAT, \texttt{PySAT}, see \cite{Ignatiev_PySAT_2018}.
To mitigate the exponential search space, we applied static symmetry breaking using \texttt{BreakID}, see the account \cite{DevriendtEtAl_BreakID_2016}, to prune isomorphic configurations under the action of the relevant
dihedral groups ($D_8$ for the square lattice) combined with global color inversion. This approach allows us to resolve the $\Xi_{14}$ case for squares and to explore significantly larger domains for some special rectangles, hexagons and cubes than were accessible to the managed brute force method. This led us to the confirmation of Walton and Li's result, \cite{WL}:

\begin{theo}[Walton - Li]\label{th:Xi}
The formulas $\Xi_m$ are satisfiable if and only if $m \leq 14$. Consequently $\Gamma(\square,2)=15$. One cannot color $S_{15}$ without the occurrence of a horizontal monochromatic square. This is the best bound for the homothety class of horizontal squares. 
\end{theo}

\begin{figure}[ht]
  \centering
  \begin{tikzpicture}[scale=0.35]
    \foreach \x/\y in {
      0/12, 0/9, 0/8, 0/7, 0/6, 0/3, 0/1,
      1/10, 1/9, 1/7, 1/5, 1/3, 1/2,
      2/12, 2/11, 2/10, 2/8, 2/7, 2/2, 2/1,
      3/13, 3/12, 3/7, 3/6, 3/4, 3/2, 3/0,
      4/12, 4/11, 4/9, 4/7, 4/5, 4/4,
      5/12, 5/10, 5/9, 5/4, 5/3, 5/1, 5/0,
      6/9, 6/8, 6/6, 6/5, 6/4, 6/2, 6/1,
      7/13, 7/11, 7/9, 7/7, 7/6, 7/1, 7/0,
      8/12, 8/11, 8/6, 8/5, 8/3, 8/2, 8/1,
      9/11, 9/10, 9/8, 9/7, 9/6, 9/4, 9/3,
      10/13, 10/12, 10/11, 10/9, 10/8, 10/3, 10/2, 10/0,
      11/13, 11/8, 11/7, 11/5, 11/4, 11/3, 11/1, 11/0,
      12/13, 12/12, 12/10, 12/8, 12/6, 12/5, 12/0,
      13/13, 13/11, 13/10, 13/5, 13/4, 13/2, 13/0}
      \fill (\x,\y) circle (0.18);

    \foreach \x/\y in {
      0/13, 0/11, 0/10, 0/5, 0/4, 0/2, 0/0,
      1/13, 1/12, 1/11, 1/8, 1/6, 1/4, 1/1, 1/0,
      2/13, 2/9, 2/6, 2/5, 2/4, 2/3, 2/0,
      3/11, 3/10, 3/9, 3/8, 3/5, 3/3, 3/1,
      4/13, 4/10, 4/8, 4/6, 4/3, 4/2, 4/1, 4/0,
      5/13, 5/11, 5/8, 5/7, 5/6, 5/5, 5/2,
      6/13, 6/12, 6/11, 6/10, 6/7, 6/3, 6/0,
      7/12, 7/10, 7/8, 7/5, 7/4, 7/3, 7/2,
      8/13, 8/10, 8/9, 8/8, 8/7, 8/4, 8/0,
      9/13, 9/12, 9/9, 9/5, 9/2, 9/1, 9/0,
      10/10, 10/7, 10/6, 10/5, 10/4, 10/1,
      11/12, 11/11, 11/10, 11/9, 11/6, 11/2,
      12/11, 12/9, 12/7, 12/4, 12/3, 12/2, 12/1,
      13/12, 13/9, 13/8, 13/7, 13/6, 13/3, 13/1}
      \draw (\x,\y) circle (0.18);
  \end{tikzpicture}
  \caption{A 2-coloring of the $14 \times 14$ grid with no monochromatic horizontal squares. This is a maximal solution for this problem. But as a $14 \times 14$ grid can be divided in four disjoint $7 \times 7$ grids, and as $\Sigma(\square, 2)=7$, there are at least four skew monochromatic squares.}
  \label{fig:square15}
\end{figure}
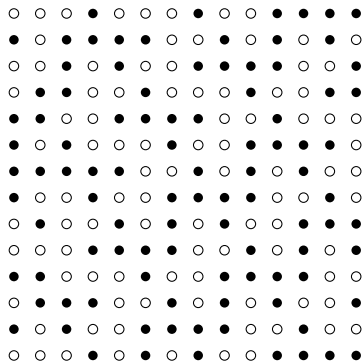

{\bf Proof}: For $m = 14$, the SAT instance corresponding to $\Xi_{14}$ is satisfiable. A corresponding $2$-coloring of the $14 \times 14$ grid with no monochromatic axes-parallel square can be constructed by the SAT-based approach, see Figure~\ref{fig:square15}. For $m = 15$, we generated the SAT instance corresponding to $\Xi_{15}$ and passed it
to the solver \texttt{CaDiCaL}, which produced a DRAT (Delete Resolution Asymmetric
Tautology) proof of unsatisfiability. The proof trace has size about $15$~GB and
certifies that no $2$-coloring of the $15\times 15$ grid avoids monochromatic
axes-parallel squares. Hence $\Xi_{15}$ is unsatisfiable.
 \qed 

%%%%%%%%%%%%%%%%%%%%%%%%%%%%%%%%%%%%%%%%%%%%%%%%%%%%%%%%%%%%%%%%%%%%%%%%%%%%%%%%%%%%%%%%%%%%%%%%%%%%%%%%%%%%%%%%%%%

\section{Rectangles} 

In this section we study rectangles in the square grids $S_m$. Fix an integer $k\ge 2$ and consider the $4$-point configuration
\[
  R_k := \{(0,0),(1,0),(0,k),(1,k)\} \subset \Z^2.
\] 

\medskip
\noindent{\bf Similarity case (all orientations).}
We want to forbid monochromatic rectangles in $S_m$ which are \emph{similar} (in the Euclidean plane) to $R_k$, in arbitrary orientations.

Let $(a,b)\in\Z^2$ with $(a,b)\neq(0,0)$ and define the two vectors
\[
u=(a,b),\qquad v=(kb,-ka).
\]
Then $u\cdot v=0$ and $\|v\|=k\|u\|$, hence for every base point $(x,y)\in\Z^2$
the four points
\[
(x,y),\quad (x,y)+u,\quad (x,y)+v,\quad (x,y)+u+v
\]
are the vertices of a rectangle similar to $R_k$.

For $m\ge 1$ we define the family of all such rectangles contained in $S_m$:
\[
  C_m^{(k)} =
  \left\{
  \{(x,y),(x+a,y+b),(x+kb,y-ka),(x+a+kb,y+b-ka)\}\subseteq S_m
  \,\middle|\,
  x,y\in\N,\ a,b\in\Z,\ (a,b)\neq(0,0)
  \right\}.
\]
In analogy with the previous sections, define
\[
  \Omega_m^{(k)} := \Upsilon(S_m, C_m^{(k)}).
\]
Thus $\Omega_m^{(k)}$ is satisfiable if and only if there exists a $2$-coloring of $S_m$
with no monochromatic rectangle similar to $R_k$.

\medskip
\noindent{\bf Homothety case (axis-parallel).}
A $\Z$-homothetic image of $R_k$ has the form
\[
  Q_k(i,j,d) := (i,j) + d\cdot R_k
  = \{(i,j),(i+d,j),(i,j+kd),(i+d,j+kd)\},
\]
where $(i,j)\in\Z^2$ and $d\in\N\setminus\{0\}$. Thus $d$ and $kd$ are the horizontal
and vertical side lengths, respectively. For $k\ge 2$ and $m\ge 1$, consider the set
\[
  A_m^{(k)} =
  \{\, Q_k(i,j,d) \mid d\in\N\setminus\{0\},\ i,j\in\N,\ i+d<m,\ j+kd<m \,\}.
\]

Also, we consider the rotated rectangle
\[
  R_k' = \{ (0,0), (k,0), (0,1), (k,1) \},
\]
and its $\Z$-homothetic images
\[
  Q_k'(i,j,d) := (i,j) + d\cdot R_k'
  = \{(i,j),(i+kd,j),(i,j+d),(i+kd,j+d)\},
\]
with horizontal side-length $kd$ and vertical side-length $d$. Their set is
\[
  B_m^{(k)} =
  \{\, Q_k'(i,j,d) \mid d\in\N\setminus\{0\},\ i,j\in\N,\ i+kd<m,\ j+d<m \,\}.
\]

In analogy with the equilateral triangles and with the squares, define
\[
  \Theta_m^{(k)} := \Upsilon(S_m, A_m^{(k)} \cup B_m^{(k)}).
\]
The following formulas $\Xi^{(k)}_m$ correspond to the $\Z$-homothetic images of the rectangle $R_k$:
\[
  \Xi_m^{(k)} := \Upsilon(S_m, A_m^{(k)}).
\]

In our implementation we generate $\Xi^{(k)}_m$, $\Theta^{(k)}_m$ and $\Omega_m^{(k)}$ in DIMACS CNF format and solve them
with CDCL SAT solvers (\texttt{Kissat} and \texttt{CaDiCaL}), using the Python library \texttt{PySAT} as an interface. For satisfiable instances we output a witness coloring. For the first unsatisfiable instance we additionally generate a DRAT certificate
of unsatisfiability, which can be independently checked by a DRAT checker.

{\bf Safe symmetry breaking}:
Since global color inversion maps any satisfying coloring to another satisfying
coloring, one may safely fix $x_{0,0}=0$ by adding the unit clause $\neg x_{0,0}$. This removes the color-flip symmetry without excluding any satisfiable grid size. 

\begin{theo}\label{thm:obl-rect-sigma}
The first unsatisfiable instances of $\Omega_m^{(k)}$ occur at
\[
m=8 \text{ for } k=2,\qquad
m=13 \text{ for } k=3,\qquad
m=15 \text{ for } k=4.
\]
Consequently,
\[
\Sigma(R_2,2)=8,\qquad \Sigma(R_3,2)=13,\qquad \Sigma(R_4,2)=15.
\]
So one cannot color these grids without occurrences of monochromatic rectangles similar to $R_k$, in arbitrary orientations.
\end{theo}

\begin{theo}\label{thm:computed-gamma-swapped}
The smallest grids for $\Theta^{(k)}_m$ to be unsatisfiable occur at:
\[
  m = 23 \,\,\text{ for } k=2,\qquad
  m = 27 \,\,\text{ for } k=3,\qquad
  m = 28 \,\,\text{ for } k=4.
\] 
So one cannot color these grids without occurrences of monochromatic axis-parallel rectangles homothetic to $R_k$ or to $R_k'$, i.e.\ in both orientations.
\end{theo} 

\begin{table}[t]
\centering
\begin{tabular}{c|c|r|r|r}
$k$ & first UNSAT $m$ & vars $=m^2$ & rectangles & clauses \\
\hline
2 & 23 & 529 & 4554 & 9109 \\
3 & 27 & 729 & 5112 & 10225 \\
4 & 28 & 784 & 4256 & 8513
\end{tabular}
\caption{Minimal UNSAT instance sizes $m$ for $\Theta^{(k)}_m$ at $k = 2, 3, 4$. The number of considered rectangles is twice as big as in $\Xi^{(k)}_m$. Clauses include the singleton clause $\neg x_{0,0}$.}
\label{tab:rect-results-swapped}
\end{table}

\begin{theo}\label{thm:computed-gamma}
The smallest grids for $\Xi^{(k)}_m$ to be unsatisfiable occur at following values of $m$:
\[
  \Gamma(R_2,2)=27,\qquad
  \Gamma(R_3,2)=40,\qquad
  \Gamma(R_4,2)=52,\qquad
  \Gamma(R_5,2)=66.
\] 
So one cannot color these grids without occurrences of monochromatic axis-parallel rectangles homothetic to $R_k$. These are the best bounds for Gallai's Theorem for the corresponding homothety classes. 
\end{theo} 

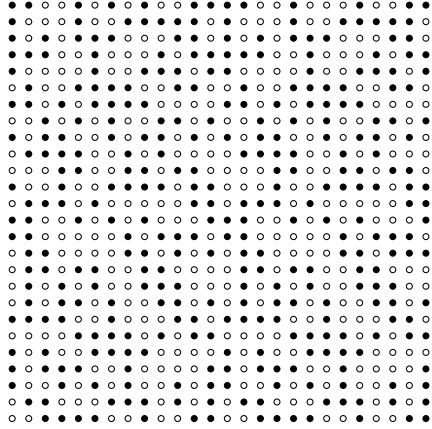
\begin{figure}[t]
  \centering
  \begin{tikzpicture}[x=0.22cm,y=0.22cm]
    \foreach \x in {0,...,25}
      \foreach \y in {0,...,25}
        \draw (\x,\y) circle (0.18);

    \foreach \x/\y in {
      2/0, 3/0, 4/0, 5/0, 8/0, 11/0, 12/0, 15/0, 16/0, 17/0, 18/0, 21/0, 24/0, 25/0,
      1/1, 2/1, 6/1, 7/1, 8/1, 10/1, 12/1, 14/1, 15/1, 19/1, 20/1, 21/1, 23/1, 25/1,
      0/2, 3/2, 5/2, 7/2, 10/2, 12/2, 13/2, 16/2, 18/2, 20/2, 23/2, 25/2,
      0/3, 2/3, 3/3, 4/3, 7/3, 12/3, 13/3, 15/3, 16/3, 17/3, 20/3, 25/3,
      0/4, 2/4, 5/4, 6/4, 7/4, 8/4, 13/4, 15/4, 18/4, 19/4, 20/4, 21/4,
      4/5, 5/5, 6/5, 7/5, 9/5, 11/5, 12/5, 17/5, 18/5, 19/5, 20/5, 22/5, 24/5, 25/5,
      0/6, 1/6, 2/6, 3/6, 6/6, 10/6, 11/6, 13/6, 14/6, 15/6, 16/6, 19/6, 23/6, 24/6,
      1/7, 3/7, 4/7, 6/7, 9/7, 10/7, 12/7, 14/7, 16/7, 17/7, 19/7, 22/7, 23/7, 25/7,
      1/8, 3/8, 5/8, 8/8, 9/8, 10/8, 14/8, 16/8, 18/8, 21/8, 22/8, 23/8,
      1/9, 2/9, 4/9, 5/9, 8/9, 9/9, 14/9, 15/9, 17/9, 18/9, 21/9, 22/9,
      1/10, 2/10, 7/10, 8/10, 10/10, 12/10, 14/10, 15/10, 20/10, 21/10, 23/10, 24/10, 25/10,
      0/11, 1/11, 7/11, 9/11, 10/11, 11/11, 13/11, 14/11, 20/11, 22/11, 23/11, 24/11,
      0/12, 1/12, 4/12, 6/12, 8/12, 12/12, 13/12, 14/12, 17/12, 19/12, 21/12, 25/12,
      1/13, 2/13, 3/13, 5/13, 11/13, 12/13, 14/13, 15/13, 16/13, 18/13, 24/13, 25/13,
      0/14, 3/14, 6/14, 7/14, 8/14, 9/14, 11/14, 12/14, 16/14, 19/14, 20/14, 21/14, 22/14, 24/14, 25/14,
      3/15, 4/15, 7/15, 8/15, 10/15, 11/15, 16/15, 17/15, 20/15, 21/15, 23/15, 24/15,
      1/16, 2/16, 3/16, 4/16, 7/16, 9/16, 14/16, 15/16, 16/16, 17/16, 20/16, 22/16,
      0/17, 2/17, 3/17, 6/17, 8/17, 9/17, 11/17, 13/17, 15/17, 16/17, 19/17, 21/17, 22/17, 24/17,
      2/18, 4/18, 6/18, 9/18, 10/18, 12/18, 15/18, 17/18, 19/18, 22/18, 23/18, 25/18,
      0/19, 1/19, 3/19, 5/19, 6/19, 7/19, 8/19, 13/19, 14/19, 16/19, 18/19, 19/19, 20/19, 21/19,
      1/20, 4/20, 5/20, 6/20, 7/20, 10/20, 11/20, 14/20, 17/20, 18/20, 19/20, 20/20, 23/20, 24/20,
      0/21, 5/21, 8/21, 9/21, 10/21, 12/21, 13/21, 18/21, 21/21, 22/21, 23/21, 25/21,
      0/22, 1/22, 2/22, 5/22, 9/22, 11/22, 12/22, 13/22, 14/22, 15/22, 18/22, 22/22, 24/22, 25/22,
      0/23, 2/23, 4/23, 5/23, 6/23, 10/23, 11/23, 13/23, 15/23, 17/23, 18/23, 19/23, 23/23, 24/23,
      0/24, 4/24, 7/24, 8/24, 9/24, 10/24, 11/24, 13/24, 17/24, 20/24, 21/24, 22/24, 23/24, 24/24,
      0/25, 1/25, 4/25, 6/25, 8/25, 11/25, 12/25, 13/25, 14/25, 17/25, 21/25, 24/25, 25/25
    } \fill (\x,\y) circle (0.18);
  \end{tikzpicture}
  \caption{A coloring of $S_{26}$ avoiding monochromatic axis-parallel rectangles of side-lengths $(d,2d)$. This is the maximal solution of the family $\Xi^{(2)}$.}
  \label{fig:k2n26}
\end{figure}

For reference, the UNSAT instances at $m=\Gamma(R_k,2)$ have $m^2$ variables and $2\cdot(\#\text{rectangles})+1$ clauses. The additional $+1$ accounts for $\neg x_{0,0}$. Figure~\ref{fig:k2n26} visualizes a maximal coloring for $k=2$ and $m=26$, i.e. a satisfying assignment for $\Xi^{(2)}_{26}$.

\begin{table}[t]
\centering
\begin{tabular}{c|c|r|r|r}
$k$ & $m=\Gamma(R_k,2)$ & vars $=m^2$ & rectangles & clauses \\
\hline
2 & 27 & 729  & 3744  & 7489 \\
3 & 40 & 1600 & 8697  & 17395 \\
4 & 52 & 2704 & 14768 & 29537 \\
5 & 66 & 4356 & 24687 & 49375
\end{tabular} 
\caption{Minimal UNSAT instance sizes at $m=\Gamma(R_k,2)$ for $\Xi^{(k)}_m$ at $k=2,3,4,5$.  Clauses include the singleton clause $\neg x_{0,0}$.}
\label{tab:rect-results-oneway}
\end{table}

%%%%%%%%%%%%%%%%%%%%%%%%%%%%%%%%%%%%%%%%%%%%%%%%%%%%%%%%%%%%%%%%%%%%%%%%%%%%%%%%%%%%%%%%%%%%%%%%%%%%%%%%%%%%%%%%%%%

\section{Symmetry classes}\label{sec:symmetry} 

In order to correctly classify the different solutions, we identify those solutions which are in the same orbit of the action of the dihedral group $D_6$ of the equilateral triangle, respectively $D_8$ of the square. 

The group $D_6$ is generated by an axis symmetry $\sigma$ and by a central rotation $\rho$ and has $6$ elements. On the triangular $m$-grid $T_m$, indexed starting with $0$, we may
take
\[
  \sigma(x,y) = (y - x, x), \qquad
  \rho(x,y) = (m - 1 - y, m - 1 - y + x).
\]
We say that two solutions $s_1(T_m)$ and $s_2(T_m)$ are $D_6$-{\bf equivalent} if for some
$\lambda \in D_6$ and for all $0 \leq x \leq  y < m$ one has
$s_1(x,y) = s_2 (\lambda(x,y))$.

The group $D_8$ is generated by an axis symmetry $\sigma$ and by a central rotation $\rho$ and has $8$ elements. On the square $m$-grid $S_m$, indexed starting with $0$, we may
take
\[
  \sigma(x,y) = (y, x), \qquad
  \rho(x,y) = (m - 1 - y, x).
\]
We say that two solutions $s_1(S_m)$ and $s_2(S_m)$ are $D_8$-{\bf equivalent} if for some $\lambda \in D_8$ and for all $0 \leq x, y < m$ one has
$s_1(x,y) = s_2 (\lambda(x,y))$.

As between $0$ and $1$ no color is privileged, we add the flip-color function
$f: \{0,1\} \rightarrow \{0,1\}$ with $f(0) = 1$ and $f(1) = 0$.  We say that two
solutions $s_1(T_m)$ and $s_2(T_m)$ are {\bf equivalent} if either they are $D_6$-equivalent
or $s_1(T_m)$ and $f \circ s_2(T_m)$ are $D_6$-equivalent. We say that two
solutions $s_1(S_m)$ and $s_2(S_m)$ are {\bf equivalent} if either they are $D_8$-equivalent
or $s_1(S_m)$ and $f \circ s_2(S_m)$ are $D_8$-equivalent.

The maximal solutions found by managed brute force have been compared under the actions of $D_6$, respectively $D_8$ and under the global color flip $f$.

\begin{theo}\label{theoremsymmetry} According to symmetry, the maximal solutions of the problems $\Phi$, $\Gamma$, $\Psi$, $\Pi$ and $\Omega$ behave as follows:
\begin{enumerate}
    \item The $18$ solutions of $\Phi_3$ fall in four $D_6$-equivalence classes: two of $6$ elements and two of $3$ elements. If we add the color flip, we get two equivalence classes: one of $12$ elements and one of $6$ elements. 
    \item The $36$ solutions of $\Gamma_4$ fall in six six-element $D_6$-equivalence classes. If we add the color flip, these equivalence classes collapse to three twelve-element equivalence classes. As the solutions for $\Gamma_4$ are identical with the solutions for $\Psi_4$, we solved the equivalence problem for these solutions as well. 
    \item The $8$ solutions of $\Pi_4$ group in two four-element $D_8$-equivalence classes. If we add the color flip, the equivalence classes do not change. 
    \item The $56$ solutions of $\Omega_6$ fall in six eight-element $D_8$-equivalence classes and in further two four-element $D_8$-equivalence classes. If we add the  color flip, the six eight-element classes collapse in three sixteen-element  classes, while the two four-element classes remain unchanged. 
\end{enumerate} 
\end{theo}

We also  investigated the structure of the solution space for $\Xi_m$ for small $m$ by counting equivalence classes under the
action of the dihedral group $D_8$ together with color flip. For $m = 4$ there are exactly $727$ distinct equivalence classes of solutions of $\Xi_4$. For
$m = 5$ we found exactly $48\,974$ distinct classes. For $m = 6$ the number of classes already exceeds $500\,000$ before the computation was halted. This illustrates a severe
combinatorial explosion in the solution space, in contrast with the abrupt phase transition to unsatisfiability at $m = 15$.

%%%%%%%%%%%%%%%%%%%%%%%%%%%%%%%%%%%%%%%%%%%%%%%%%%%%%%%%%%%%%%%%%%%%%%%%%%%%%%%%%%%%%%%%%%%%%%%%%%%%%%%%%%%%%%%%%%%

\section{Lower bounds for hexagons and cubes}\label{sec:hex-cubes}

 We first extend the investigation to regular hexagons on the triangular lattice. We map the lattice points to a coordinate system $\Z^2$ where valid points $(x,y)$ satisfy
$x \equiv y \pmod 2$. In this geometry, a  regular hexagon of side length $s$ is defined as the set of six vertices
\[
  H(x,y,s) = \{ (x,y), (x+2s, y), (x, y+2s), (x+2s, y+2s), (x-s, y+s), (x+3s, y+s) \},
\]
expressed in suitable coordinates.  In our SAT encoding for hexagons, the parameter $m$ is the
size of a rectangular window
$R_m^* = \{(x,y) \in \mathbb{Z}^2 : 0 \le x < 2m,\ 0 \le y < m,\ x \equiv y \pmod 2\}$. 

The set of interest is
$$\hexagon_m =\{ H(x,y,s) \,|\, s \in \N, s > 0, H(x,y,s) \subseteq R_m^* \},$$
and the propositional formula expressing that no such hexagon is monochromatic is given by:
$$H_m = \Upsilon(R_m^*, \hexagon_m).$$

We define the {\bf Gallai Number} for regular hexagons to be the minimal number $m$ such that no matter how the elements of $R_m^*$ are $2$-colored, a monochromatic axes-parallel hexagon arises.

\begin{theo}\label{th:hex} The formula $H_{94}$ is satisfiable. The Gallai Number for regular hexagons on the triangular lattice is at least $95$.
\end{theo}

Each row has $m$ valid points, so $|R_m^*| = m^2$; in particular $|R_{94}^*| = 94^2 = 8836$.
A regular hexagon of side length $s$ has vertical span $2s$, so we must have $s \le \lfloor (m-1)/2 \rfloor$,
hence its diameter is at most $2\lfloor (m-1)/2 \rfloor < m$.
For $m=94$ this gives $s \le 46$ and diameter $\le 92$.
Thus $H_{94}$ satisfiable means there is a $2$-coloring of $R_{94}^*$ with no monochromatic hexagon of diameter up to $92$.

This stands in marked contrast to the square lattice case, where the corresponding
Gallai number for axes-parallel squares is $15$. The triangular lattice appears to offer substantially more combinatorial freedom to avoid monochromatic structures,
despite its higher local connectivity, because a regular hexagon involves six vertices
instead of four.

We next consider the $3$-dimensional analogue. We search for monochromatic
axes-parallel cubes in the discrete grid $\{0,\dots,n-1\}^3$. A cube of side length $s$
is determined by $8$ vertices, given by the Cartesian product of the sets $\{x, x+s\}$, $\{y, y+s\}$, and $\{z, z+s\}$. 

\begin{theo}\label{th:cubes} For three-dimensional cubes in $\Z^3$, $\Gamma(\{0,1\}^3, 2) \geq 181$.
\end{theo}

{\bf Proof}: Applying the SAT encoding with symmetry breaking as above, we found a
$2$-coloring of the $180 \times 180 \times 180$ grid with no monochromatic axes-parallel cube. This shows that the Gallai number for cubes is at least $181$. At
the time of writing we do not yet have an upper bound, so the exact value remains open.
\qed

The large lower bound is consistent with a simple probabilistic consideration: a cube consists of $8$ vertices, so in a random $2$-coloring of the infinite lattice the probability that a given cube is monochromatic is only $2^{-7} = 1/128$. This relatively weak local constraint allows for very large grids to be colored without forcing a monochromatic cube.

%%%%%%%%%%%%%%%%%%%%%%%%%%%%%%%%%%%%%%%%%%%%%%%%%%%%%%%%%%%%%%%%%%%%%%%%%%%%%%%%%%%%%%%%%%%%%%%%%%%%%%%%%%%%%%%%%%%%%%%%

\section{Implications in discrete geometry}

Here are some consequences for colorings of some
infinite lattices and for colorings of the Euclidean plane.

\begin{cor}
For every coloring $s : T \rightarrow \{0,1\}$ of the infinite triangular lattice $T$, the colored lattice contains:
\begin{enumerate}
\item Infinitely many monochromatic equilateral triangles of edge-length $\leq 3$, occurring in various orientations.
\item Infinitely many monochromatic equilateral triangles, which are parallel with the
      lattice axes and are upwards oriented, of edge-length $\leq 4$.
\end{enumerate}
Moreover, there are infinitely many larger monochromatic equilateral triangles of each kind, and their edge lengths have no upper bound.
\end{cor}

{\bf Proof}: The first claim is true because the formulas $\Phi_m$ are satisfiable if
and only if $m \leq 3$. This means that, no matter how one colors $T_4$, there will be
a monochromatic equilateral triangle with vertices in its elements. But the edge-length
of $T_4$ is $3$. The second claim comes from the fact that $\Psi_m$ is satisfiable if
and only if $m \leq 4$. For the fact that there is no upper bound for the edge-length of monochromatic triangles, observe that for all $k \in \N$, we can consider instead the triangular lattice $k \Z + k \Z \theta$.  \qed

\begin{cor}
For every coloring $s : \Z \times \Z \rightarrow \{0,1\}$ of the square lattice, the
lattice contains:
\begin{enumerate}
\item Infinitely many squares, which are parallel with the lattice axes and have an even number of black vertices, of edge-length $\leq 2$.
\item Infinitely many squares, which are parallel with the lattice axes and have a different number of black and white vertices, of edge-length $\leq 4$.
\item Infinitely many monochromatic horizontal or skew squares, of edge-length $\leq 6$.
\item Infinitely many monochromatic horizontal squares of edge-length $\leq 14$.
\end{enumerate}
Moreover, there are infinitely many larger squares of each kind, and their 
edge-lengths have no upper bound.
\end{cor}

{\bf Proof}: The first claim is a consequence of the linear problem considered as appetizer. The second one follows from Theorem~\ref{th:Pi} concerning the formulas
$\Pi_m$ and the third one from Theorem~\ref{th:Omega} concerning the formulas $\Omega_m$. By Theorem~\ref{th:Xi}, every $15\times 15$ subgrid of $\Z^2$ contains a
monochromatic axes-parallel square. The infinite lattice $\Z^2$ contains infinitely many
pairwise disjoint $15\times 15$ sub-grids, so we obtain infinitely many such squares.
Their edge length is at most $14$, since this is the maximal side length of a square
contained in a $15\times 15$ grid.

Again, the square lattice contains subgroups $k\Z \times k\Z$ with arbitrarily large unit distance, and those results apply for these lattices as well. \qed 

The following statement, similar with statements proved by Alm in \cite{A}, is a consequence of Gallai's Theorem. In the cited work, Alm did not observe the density aspect.

\begin{cor}
    If $U \subset \R^2$ is a simply connected open set, and $s : \R^2 \rightarrow \{0, 1, \dots, k-1\}$ is some $k$-coloring of the plane, then $U$ contains continuum many monochromatic equilateral triangles and continuum many monochromatic squares. 
\end{cor}

{\bf Proof}: Take some equilateral triangle $T \subset U$. By Gallai's Theorem, we know that there is a minimal grid $T_m$ which cannot be $k$-colored without the occurrence of a monochromatic equilateral triangle. Divide the edges of $T$ in $m-1$ equal segments and construct a copy of $T_m$ inside $T$. Now the finite grid is $k$-colored and a monochromatic equilateral triangle occurs. Continuum many infinitesimal translations of this grid keep it inside $U$, and for any translation one has a monochromatic equilateral triangle. The proof for squares is identical. \qed 

But without any information about the minimal grid, one has no lower bound for the edge-length of these monochromatic figures. This can be improved as follows:

\begin{cor}
Let $s : \R^2 \rightarrow \{0,1\}$ be a coloring of the plane, and let $U \subset \R^2$ be a simply connected open set. 
\begin{enumerate}
\item For every equilateral triangle $T$ of edge-length $d$, ($d\in \R$, $d>0$) contained in $U$, there are continuum many monochromatic equilateral triangles contained in $U$, and every of their edge-lengths $l$ satisfies $ d \geq l \geq d/4$.  All these equilateral triangles have edges parallel with the edges of $T$ and have the same orientation as $T$. Also, there are continuum many monochromatic equilateral triangles contained in $U$, and every of their edge-lengths $l$ satisfies $ d \geq l \geq d/3$.
\item  For every square $S$ of edge-length $d$, ($d\in \R$, $d>0$) contained in $U$, there are continuum many monochromatic squares contained in $U$, and every of their edge-lengths $l$ satisfies $d \geq l\geq d/14$.  All these squares have edges parallel with the edges of $S$. Also, there are continuum many monochromatic squares contained in $U$,  and every of their edge-lengths $l$ satisfies  $d \geq l \geq d/6$. 
\end{enumerate}
\end{cor} 

{\bf Proof}: Match the grid $T_5$ on the triangle $T$ and you get a monochromatic equilateral triangle, parallel with $T$ and of the same orientation, of edge-length $\geq d/4$. The same is true for every  translation $T+\vec{x}$ such that $ T + \vec{x} \subset U$, so we get continuum many monochromatic equilateral triangles of edge-length $l$ with  $d \geq l \geq d/4$ parallel with $T$ and oriented the same way as $T$. The other statements have analogous proofs. 
\qed 

One can state similar corollaries for homothetic images of the rectangles $R_2$, $R_3$, $R_4$ or $R_5$. 

%%%%%%%%%%%%%%%%%%%%%%%%%%%%%%%%%%%%%%%%%%%%%%%%%%%%%%%%%%%%%%%%%%%%%%%%%%%%%%%%%%%%%%%%%%%%%%%%%%%%%%%%%%%%%%%%%%%

\end{document}